\documentclass[11pt,reqno]{amsart}
\usepackage{amssymb,graphicx}

\usepackage{amsfonts,amscd,amsthm,amsmath}

\usepackage{mathrsfs}

\usepackage{array}

\newcolumntype{C}[1]{>{\centering\arraybackslash}p{#1}}

\numberwithin{equation}{section}

\begin{document}

\begin{center}

$\;$\\[-8mm]\textbf{\large{Higher-order tangent and secant numbers}}

\vskip 2mm\textbf{Djurdje Cvijovi\'{c}}

\vskip 2mm {\it Atomic Physics Laboratory, Vin\v{c}a Institute of
Nuclear Sciences \\
P.O. Box $522,$ $11001$ Belgrade$,$ Republic of Serbia}\\
\textbf{E-Mail: djurdje@vinca.rs}\\

\vskip 2mm \begin{quotation} \textbf{Abstract.} In this paper higher-order tangent numbers and higher-order secant numbers, $\left\{\mathscr{T}(n,k)\right\}_{n,k\,=0}^{\infty}$ and $\left\{\mathscr{S}(n,k)\right\}_{n,k\,=0}^{\infty}$, have been studied in detail. Several known results regarding $\mathscr{T}(n,k)$ and  $\mathscr{S}(n,k)$ have been brought together along with many new results and insights and they all have been proved in a simple and unified manner. In particular, it is shown that the higher-order tangent numbers $\mathscr{T}(n,k)$ constitute a special class  of the partial multivariate Bell polynomials and  that $\mathscr{S}(n,k)$ can be  computed from the knowledge of $\mathscr{T}(n,k)$. In addition, a simple explicit formula involving a double finite sum is deduced for the numbers $\mathscr{T}(n,k)$  and it is shown that  $\mathscr{T}(n,k)$ are linear combinations of the classical tangent numbers $T_n$.
\end{quotation}
\end{center}

\vskip 1mm\noindent\textbf{2000 \textit{Mathematics Subject
Classification.}} Primary 11B83, 11B99; Secondary 11B37, 33B10.
\vskip 1mm\noindent \textbf{\textit{Key Words and Phrases.}} \begin{Small} Tangent numbers; Tangent numbers of order $k$; Higher-order (or, generalized) tangent numbers; Secant numbers; Secant numbers of order $k$;  Higher-order (or, generalized) secant numbers; Partial Bell polynomials; Stirling numbers of the second kind; Lah numbers; Central factorial numbers of the second kind; Derivative polynomials. \end{Small}

\section{Introduction}

Let us, for non-negative integers $n$ and $k$, define the {\em n--th tangent number of order $k$},  $\mathscr{T}(n,k),$ by the generating relation (see \cite[p. 259]{Comtet}; {\em cf.} \cite[p. 428]{Carlitz} and \cite[p. 305]{Carlitz2})
\begin{equation}
\frac{\tan^{k}t}{k!} = \sum_{n\, = k}^{\infty}
\mathscr{T} (n,k)\,\frac{t^n}{n!}.
\end{equation}
\noindent Similarly, define {\em n--th  secant number of order $k$}, $\mathscr{S}(n,k),$  by ({\em cf.} \cite[p. 428]{Carlitz})
\begin{equation}
\frac{\sec t \tan^{k}t}{k!}=\sum_{n\, = k}^{\infty}
\mathscr{S}(n,k)\,\frac{t^n}{n!}.
\end{equation}
\noindent Since  $\mathscr{T}(n,k)$ and  $\mathscr{S}(n,k)$ are,  respectively, generalizations of the classical (and well-known) tangent (or, zag) numbers $T_n$ (see \cite[p. 259]{Comtet} and \cite{Knuth})
\begin{equation}
\tan t = \sum_{n\, = 1}^{\infty}
T_n\,\frac{t^n}{n!},\qquad \textup{i.e.,}\;\; T_n := \mathscr{T} (n,1)
\end{equation}
\noindent and  secant (or, zig) numbers $S_n$ (see \cite{Knuth} and \cite[p. 63]{Srivastava})
\begin{equation}
\sec t = \sum_{n\, = 0}^{\infty}
S_n\,\frac{t^n}{n!},\qquad \textup{i.e.,}\;\;  S_n := \mathscr{S} (n,0),
\end{equation}
\noindent the numbers $\mathscr{T}(n,k)$ and  $\mathscr{S}(n,k)$ may also be referred to as the higher-order (or, generalized) tangent numbers and the higher-order (or, generalized) secant numbers.

In a recent paper by Cvijovi\'{c} \cite{Cvijovic}, very simple and compact closed-form higher derivative formulae are derived for eight trigonometric and hyperbolic functions   \cite[Corollaries 1 and 2]{Cvijovic} and they involve the numbers $T(n,k)$ and $S(n,k)$ given by
\begin{equation} T(n,k) = \mathscr{T}(n,k) \,k! \quad\textup{and}\quad S(n,k) = \mathscr{S}(n,k)\,k!.
\end{equation}
\noindent The elegance and remarkable simplicity of the obtained results (see Propositions 11 and 12) have been the main motivation behind our further interest in these numbers. It has turned out that they are not sufficiently studied as yet and it is aimed here to thoroughly examine $\mathscr{T}(n,k)$ and  $\mathscr{S}(n,k)$. Several known results have been brought together along with many new results and insights and they all have been proved in a simple and unified manner (for further details, see Section 3).

\section{Properties of numbers $\mathscr{T}(n,k)$ and  $\mathscr{S}(n,k)$}

In what follows, it is assumed, unless otherwise indicated, that $k, l, m$ and $n$ are non-negative integers and we set an empty sum to be zero and $$D_{t}^{n}:= \frac{d^n}{d t^n}\qquad\qquad(n\geq1).$$

We begin by an observation that it can be  easily seen, from definitions (1.1) and (1.2) and after some parity considerations, that  Proposition 1 holds  true.

\vskip 2mm \noindent{\bf Proposition 1.} {\em Let  $\mathscr{T}(n,k)$  and  $\mathscr{S}(n,k)$ be the numbers defined by $\textup{(1.1)}$ and $\textup{(1.2)}$. Then:} {\bf \textup{a)}} {\em $\mathscr{T}(n,k)\neq 0$ is only when  $1\leq k\leq n $ and either both $n$ and $k$ are even or both $n$ and $k$ are odd. In other words, $\mathscr{T}(2 m, 2 l + 1) = 0$ and $\mathscr{T}(2 m + 1, 2 l) = 0$;} {\bf \textup{b)}} {\em $\mathscr{S}(n,k)\neq 0$ is only when  $0\leq k\leq n $ and either both $n$ and $k$ are even or both $n$ and $k$ are odd. In other words, $\mathscr{S}(2 m, 2 l + 1) = 0$ and $\mathscr{S}(2 m + 1, 2 l) = 0$.}

\vskip 2mm \noindent{\bf Proposition 2.} {\em Assume that $n, k\geq1$ and let $\mathscr{T}(n,k)$ and $T_n$, respectively,  be the higher-order tangent numbers and tangent numbers. Then, we have}
\begin{equation} \mathscr{T}(k,k)= T_1^{k}
\end{equation}
\noindent {\em and, for $n\neq k$},
\begin{equation} \mathscr{T}(n,k ) = \frac{n + 1}{n-k}\cdot \frac{1}{T_1} \sum_{ r\,= 1}^{n-k} \binom{n}{r}\left(\frac{k+1}{n+1} - \frac{1}{r+1}\right) \, T_{r+1}\,\mathscr{T}(n-r,k).
\end{equation}

\vskip 2 mm \noindent {\em Demonstration.} The proof is based on the following known formula for powers of series \cite[Eqs. (1.1) and (3.2)]{Gould}. For a fixed $k$, put
\begin{equation*} \left(\sum_{n\,=1}^{\infty} a_n t^n\right)^k = \sum_{n\,= k}^{\infty} b_n t^n\qquad\qquad(k\geq1),
\end{equation*}
\noindent then the coefficients $b_n$, $n\geq k,$ are given by $b_k = a_1^k$ and
\begin{equation*}b_n=\frac{1}{a_1 (n-k)} \sum_{r\,=1}^{n-k} \Big[(r+1) (k+1)-(n+1)\Big] a_{r+1} \,b_{n-r}\quad(n\geq k+1).
\end{equation*}
It is not difficult to show that this result together with the definitions of the numbers $\mathscr{T}(n,k)$ and $T_n$ in (1.1) and (1.3) yields (2.1) and (2.2).

\vskip 2mm \noindent{\bf Proposition 3.} {\em Let $\mathscr{T}(n,k),$ $\mathscr{S}(n,k)$ and $S_n$, respectively,  be the higher-order tangent numbers, the higher-order secant numbers and secant numbers. Then}

\begin{equation}\mathscr{S}(n ,k) = \sum_{n^* = 0}^{n}\binom{n}{n^*}\mathscr{T}(n^*,k) \,S_{n - n^*}.\end{equation}

\vskip 2 mm \noindent {\em Demonstration.} Recall that if $\sum\nolimits_{n\,=0}^{\infty} a_n$ and $\sum\nolimits_{n\,= 0}^{\infty} b_n$ are two series, then their the Cauchy product is the series $\sum\nolimits_{n\,= 0}^{\infty} c_n$ where $c_n = \sum\nolimits_{k\,= 0}^{n} a_k b_{n-k}$. Hence, the required formula (2.3) is implied by the Cauchy product of power series expansions (1.4) and (1.1) for functions involving the product $\sec t \tan^k t/k!$ in (1.2).

\vskip 2mm \noindent{\bf Proposition 4.} {\em Consider the  multivariate  (exponential) partial  Bell polynomials $B_{n,k}(x_1, x_2, \ldots, x_n)$ defined by the formal power series expansion } \textup{(} \cite[pp. 133--137]{Comtet} {\em and}  \cite[pp. 412--417]{Charalambides}\textup{)}
\begin{equation} \frac{1}{k!} \left(\sum_{m\,= 1}^{\infty} x_m \frac{t^m}{m!}\right)^k = \sum_{n\,= k}^{\infty} B_{n,k}(x_1, x_2, \ldots, x_n)\,\frac{t^n}{n!}
\end{equation}
\noindent {\em and let $\mathscr{T}(n,k)$ and $T_n$ be the higher-order tangent numbers and tangent numbers, respectively. Then}
\begin{equation}\mathscr{T}(n,k) = B_{n,k}(T_1, T_2, \ldots, T_n).
\end{equation}

\vskip 2 mm \noindent {\em Demonstration.} The proposed identity is readily available from (2.4) and (1.3) in conjunction with (1.1).

\vskip 2 mm In view of Proposition 4, the higher-order tangent numbers constitute a special class  of the Bell polynomials. Thus, for fixed $n$ and $k$, $\mathscr{T}(n,k)$ is a homogenous  and isobaric polynomial in $T_1, T_2, \ldots, T_n$ of total degree $k$ and total weight $n$, i.e., it is a linear combination of monomials $T_1^{k_1}\cdot T_2^{k_2}\cdots T_n^{k_n}$ whose partial degrees and weights are constantly given by $k_1+k_2+\ldots + k_n =k$ and $k_1 + 2 k_2+\ldots + n k_n = n$.

\vskip 2mm \noindent{\bf Example 1.} In order to demonstrate an application of Proposition 4, we tabulate several  $\mathscr{T}(n,k)$ given in terms of $T_n$. Note that a list of $B_{n,k} \equiv B_{n,k}(x_1, x_2,\ldots, x_n)$ for $k\leq n\leq12$ can be, for instance, found in \cite[pp. 307--308]{Comtet}, while the needed values of $T_n$ are $T_1 = 1, T_3 = 2, T_5 = 16, T_7 = 272, T_9 = 7936.$
\begin{small}
\begin{align*}&\mathscr{T}(1,1) = T_1,
\\
&\mathscr{T}(2,2) = T_1^2,
\\
&\mathscr{T}(3,1) = T_3,\quad\mathscr{T}(3,3) = T_1^3,
\\
&\mathscr{T}(4,2) = 4\, T_1 T_3,\quad\mathscr{T}(4,4) = T_1^4,
\\
&\mathscr{T}(5,1) = T_5,\quad\mathscr{T}(5,3) =  10\, T_1^2 T_3,\quad\mathscr{T}(5,5) = T_1^5,
\\
&\mathscr{T}(6,2) = 6\, T_1 T_5 + 10\, T_3^2,\quad\mathscr{T}(6,4) = 20\, T_1^3 T_3,\quad\mathscr{T}(6,6) = T_1^6,
\\
&\mathscr{T}(7,1) = T_7,\quad\mathscr{T}(7,3) =  21\, T_1^2 T_5 + 70\, T_1 T_3^2,\quad\mathscr{T}(7,5) = 35\, T_1^4 T_3,\quad\mathscr{T}(7,7) = T_1^7,
\\
&\mathscr{T}(8,2) =  8\, T_1 T_7 + 56\, T_3 T_5,\quad\mathscr{T}(8,4) = 56\, T_1^3 T_5 + 280\, T_1^2 T_3^2,\quad\mathscr{T}(8,6) = 56\, T_1^5 T_3,
\\
&\hskip 15mm \mathscr{T}(8,8) = T_1^8,
\\
&\mathscr{T}(9,1) = T_9,\,\mathscr{T}(9,3) =  36\, T_1^2 T_7 + 504\, T_1 T_3 T_5 + 280\, T_3^3,\,\mathscr{T}(9,5) = 126\, T_1^4 T_5 + 840\, T_1^3 T_3^3,
\\
&\hskip 15mm\mathscr{T}(9,7) = 84 T_1^6 T_3,\, \mathscr{T}(9,9) = T_1^9.
\end{align*}
\end{small}

\noindent{\bf Corollary 1.} {\em Let $\mathscr{T}(n,k)$ and $T_n$ be the higher-order tangent numbers and tangent numbers, respectively. Then}
\begin{align}
&\textup{a)}\quad  \mathscr{T}(n+1 ,k+1) = \sum_{r \,=  0}^{n-k}\binom{n}{r} T_{r+1}\,\mathscr{T}(n-r,k),\hskip 30mm
\\
&\qquad\mathscr{T}(0,0) = 1,\qquad\mathscr{T}(n,0) = 0\qquad(n\geq 1);\nonumber
\\[2mm]
&\textup{b)}\quad  \mathscr{T}(n+1 ,k+1) = \frac{1}{k+1 }\sum_{r \,=  0}^{n-k}\binom{n+1}{r+1} T_{r+1}\,\mathscr{T}(n-r,k),
\\
&\qquad\mathscr{T}(0,0) = 1,\qquad\mathscr{T}(n,0) = 0\qquad(n\geq 1).\nonumber
\end{align}

\vskip 2 mm \noindent {\em Demonstration.} The results are immediate consequences of Proposition 4 and the recurrence relations \cite[p. 415, Eqs. (11.11) and (11.12)]{Charalambides}
\begin{equation*}B_{n + 1,k + 1} = \sum_{r\, = 0}^{n-k} \binom{n}{r}  x_{r+1}\,B_{n-r,k}
\end{equation*}
\noindent and
\begin{equation*}B_{n + 1,k + 1} = \frac{1}{k+1} \sum_{r\, = 0}^{n-k} \binom{n+1}{r+1}  x_{r+1}\,B_{n-r,k}.
\end{equation*}
\noindent satisfied by  the partial Bell polynomials $B_{n,k}$ with $B_{0,0}:=1$.

\vskip 2mm \noindent{\bf Proposition 5.} {\em Let  $S_n$ be the secant numbers $\textup{(1.4)}$. The  numbers $\mathscr{T}(n,k)$  and $\mathscr{S}(n,k)$ obey the following recurrence relations}
\begin{equation} \textup{a)}\quad \mathscr{T}(n + 1 ,k) = \mathscr{T}(n,k - 1) + k \,(k+1) \mathscr{T}(n,k+1)\quad(n\geq0, k\geq1),\end{equation}
$$\mathscr{T}(0,0) = 1,\qquad\mathscr{T}(n,0) = 0\qquad(n\geq 1);$$
\begin{equation}\textup{b)}\quad \mathscr{S}(n + 1 ,k) = \mathscr{S}(n,k - 1) +  (k + 1)^2
\mathscr{S}(n,k+1)\quad(n\geq0, k\geq1),
\end{equation}
$$\mathscr{S}(n,0)= S_n\qquad(n\geq0).$$

\vskip 2 mm \noindent {\em Demonstration.} The recurrence relation (2.8) follows at once from
\begin{equation*} \frac{\tan^{k-1}t}{(k-1)!} + k (k+1) \frac{\tan^{k + 1}t}{(k +1)!} = \frac{  \sec^2 t \tan^{k-1} t}{(k-1)!} =  \frac{1}{k!} \, D_{t} \tan^k t
\end{equation*}
\noindent and
\begin{equation*}
\frac{1}{k!} \, D_{t} \tan^k t = \sum_{n+1\,= k}^{\infty} \mathscr{T}(n + 1,k) \frac{t^n}{n!}.
\end{equation*}
\noindent Similarly, by making use of
\begin{equation*}
\frac{1}{k!} \, D_{t} \big( \sec t \tan^k t \big) = \frac{  \sec^3 t \tan^{k-1} t}{(k-1)!} + \frac{  \sec t \tan^{k + 1} t}{k!} = \sum_{n+1\,=k}^{\infty} \mathscr{T}(n + 1,k) \frac{t^n}{n!}.
\end{equation*}
\noindent and
\begin{equation*} \frac{\sec t \tan^{k-1}t}{(k-1)!} + k (k+1) \frac{\sec t \tan^{k + 1}t}{(k +1)!} = \frac{  \sec^3 t \tan^{k-1} t}{(k-1)!},
\end{equation*}
\noindent it is straightforward to arrive at (2.9).

\vskip - 2 mm \noindent \begin{table}[ht]
\caption{Higher-order tangent numbers  $\mathscr{T}(n,k)$} 
\centering 
\small
\begin{tabular}{c| C{.75cm} C{.75cm} C{.75cm} C{.75cm} C{.75cm} C{.75cm} C{.75cm} C{.75cm} C{.75cm} C{.75cm}}
$n\backslash k $ & 0 & 1 & 2 & 3 & 4 & 5 & 6 & 7 & 8 & 9 \\ [.5ex] 
\hline \\ 
0 & 1\\
1 & 0 & 1\\
2 & 0 & 0 & 1\\
3 & 0 & 2 & 0 & 1\\
4 & 0 & 0 & 8 & 0 & 1\\
5 & 0 & 16 & 0 & 20 & 0 & 1\\
6 & 0 & 0 & 136 & 0 & 40 & 0 & 1\\
7 & 0 & 272 & 0 & 616 & 0 & 70 & 0 & 1\\
8 & 0 & 0 & 3968 & 0 & 2016 & 0 & 112 & 0 & 1\\
9 & 0 & 7936 & 0 & 28160 & 0 & 5376 & 0 & 168 & 0 & 1 \\[1ex] 
\hline 
\end{tabular}
\label{table:Tab1} 
\end{table}

\vskip 2mm  For the sake of ready reference,  by employing  Proposition 5, we compute and list several of the numbers $\mathscr{T}(n,k)$ in Table 1  and  several of the numbers $\mathscr{S}(n,k)$ in Table 2.

\vskip 1 mm \noindent  \begin{table}[ht]
\caption{Higher-order secant numbers  $\mathscr{S}(n,k)$} 
\centering 
\small
\begin{tabular}{c| C{.75cm} C{.75cm} C{.75cm} C{.75cm} C{.75cm} C{.75cm} C{.75cm} C{.75cm} C{.75cm} C{.75cm}}
$n\backslash k $ & 0 & 1 & 2 & 3 & 4 & 5 & 6 & 7 & 8 & 9 \\ [.5ex] 
\hline \\ 
0 & 1\\
1 & 0 & 1\\
2 & 1 & 0 & 1\\
3 & 0 & 5 & 0 & 1\\
4 & 5 & 0 & 14 & 0 & 1\\
5 & 0 & 61 & 0 & 30 & 0 & 1\\
6 & 61 & 0 & 331 & 0 & 55 & 0 & 1\\
7 & 0 & 1385 & 0 & 1211 & 0 & 91 & 0 & 1\\
8 & 1385 & 0 & 12284 & 0 & 3486 & 0 & 140 & 0 & 1\\
9 & 0 & 50521 & 0 & 68060 & 0 & 8526 & 0 & 204 & 0 & 1\\[1ex] 
\hline 
\end{tabular}
\label{table:Tab2} 
\end{table}

\vskip 2mm \noindent{\bf Proposition 6.} {\em We have that:}

\begin{align}&\textup{a)}\quad  \sum_{k\,= 1}^{n+1} (k-1)! \mathscr{T}(n +1 ,k) = \begin{cases} \big(2^{n} -1\big) T_n \quad n \,\,\textup{odd}\\ \quad\quad\;\; 2^n S_n \quad n \,\,\textup{even} \end{cases}; \hskip 20mm
\\[5mm]
&\textup{b)}\quad \sum_{k\,= 0}^n k! \, \mathscr{T}(n,k) = 2^{n-1} \begin{cases} T_n \quad n \,\,\textup{odd}\\ S_n \quad n \,\,\textup{even} \end{cases}.
\end{align}

\newpage \noindent {\em Demonstration.} $\,$  Observe that $n$th derivative of $\tan z\;$ is a polynomial in $\tan z$, {\em i.e.}  $D_{t}^n \tan z = P_n (\tan z),$ where $P_n(x)$ is explicitly given by $\,$ ({\em cf.} (1.3), (1.5)  and \cite[Eq. (3.5)]{Cvijovic})
\begin{equation*}P_n(x) = T_n + \sum_{k\,=1}^{n+1} (k-1)! \, \mathscr{T}(n + 1,k)\, x^k.
\end{equation*}
Then, having in mind the identities $D_{t}^n \tan (\pi/4+z/2) = 2^{-n} P_n\big(\tan (\pi/4+z/2)\big)$ and $\tan (\pi/4+z/2) = \tan z + \sec z$, by the last formula  we obtain (2.10).  Further, starting from the left-hand side of (2.10) and by applying the recurrence relation (2.8), we have
\begin{equation*}\sum_{k\,= 1}^{n+1} (k-1)! \mathscr{T}(n +1 ,k) = \sum_{k\,= 0}^{n} k! \mathscr{T}(n ,k) + \sum_{k\,= 2}^{n+2} k! \mathscr{T}(n,k)= 2 \sum_{k\,= 0}^{n} k! \mathscr{T}(n ,k) - T_n,
\end{equation*}
\noindent which, by appealing to the right-hand side of (2.10), yields (2.11).

In  Proposition 7 below we shall require the {\em higher-order arctangent numbers} $\mathscr{T}^{*} (n,k)$ defined by the generating relation \cite[p. 260]{Comtet}
\begin{equation}
\frac{\arctan^{k}t}{k!} = \sum_{n\, = k}^{\infty}
\mathscr{T}^{*} (n,k)\,\frac{t^n}{n!},
\end{equation}
\noindent or, equivalently, by the recurrence relation \cite[p. 260]{Comtet}
\begin{equation}\mathscr{T}^{*} (n + 1,k) = \mathscr{T}^{*} (n,k-1) - n (n-1) \mathscr{T}^{*} (n - 1,k)\quad(n\geq0, k\geq1),
\end{equation}
$$\mathscr{T}^{*}(0,0) = 1,\qquad\mathscr{T}^{*}(n,0) = 0\qquad(n\geq 1).$$

\vskip -2 mm \noindent  \begin{table}[ht]
\caption{Higher-order arctangent numbers  $\mathscr{T}^{*}(n,k)$} 
\centering 
\small
\begin{tabular}{c| C{.2cm} C{1.2cm} C{1.2cm} C{1.3cm} C{.95cm} C{.95cm} C{.9cm} C{.9cm} C{.9cm} C{.9cm}}
$ n \backslash k $ & 0 & 1 & 2 & 3 & 4 & 5 & 6 & 7 & 8 & 9 \\ [.5ex] 
\hline \\ 
0 & 1\\
1 & 0 & 1\\
2 & 0 & 0 & 1\\
3 & 0 & --2\;\; & 0 & 1\\
4 & 0 & 0 & --8\;\; & 0 & 1\\
5 & 0 & 24 & 0 & --20\;\;& 0& 1\\
6 & 0 & 0  & 184& 0 & --40\;\;& 0& 1\\
7 & 0 &--720\;\;& 0& 784& 0& --70\;\;& 0& 1\\
8 & 0 & 0& --8448\;\;& 0& 2464& 0& --112\;\;& 0& 1\\
9 & 0 &40320& 0& --52352\;\;& 0& 6384& 0& --168\;\;& 0& 1\\[1ex] 
\hline 
\end{tabular}
\end{table}

The numbers $\mathscr{T}^{*} (n,k)$ (see Table 3) are involved in the following formula  which can be easily proved ({\em cf.} \cite[ Eqs. (3.4), (3.6) and (3.7)]{Lomont})

\begin{equation} \tan^n t = \frac{1}{(n-1)!} \sum_{r\,= 0}^{n-1} \mathscr{T}^{*} (n,r +1)\, D_{t}^r \tan t + \begin{cases} 0 \qquad\quad\, n \,\,\textup{odd} \\ (-1)^{\frac{n}{2}} \quad n \,\,\textup{even} \end{cases}\quad(n\geq1).\end{equation}

Indeed, let us make use of the induction on $n$. In case $n = 1$ the formula (2.14) is clearly true. For the induction step, we assume that (2.14) is true for $n$ and establish it for $n + 1$. This obviously follows, since by the  recurrence relation (2.13) we have
\begin{align*}\frac{1}{n!} \sum_{r \, = 0}^{n} \mathscr{T}^{*} (n  + 1, r + 1)\, D_{t}^r \tan t = \Delta_1-\Delta_2 = \frac{1}{n} \, D_{t} \tan^n t - \tan^{n-1}t = \tan^{n+1}t,
\end{align*}
\noindent where
\begin{align*}n \,\Delta_1 &= \frac{1}{(n -1) !} \sum_{r \,= 1}^{n}\mathscr{T}^{*} (n,r) D_{t}^r \tan t =  D_{t} \left( \frac{1}{(n-1)!} \sum_{r \,= 0}^{n-1}\mathscr{T}^{*} (n,r +1) D_{t}^r \tan t\right)
\\
& = D_{t} \tan^n t
\end{align*}
\noindent and
\begin{align*}(n-2)! \, \Delta_2 & = \sum_{r \,= 0}^{n}\mathscr{T}^{*} (n -1,r +1) D_{t}^r \tan t = \sum_{r \,= 0}^{n-2}\mathscr{T}^{*} (n -1,r +1)\, D_{t}^r \tan t
\\
& = (n-2)!  \tan^{n-1} t.
\end{align*}

\vskip 2mm \noindent{\bf Proposition 7.} {\em Let $T_n$ and $\mathscr{T}^{*}(n,k)$ be the tangent numbers and the higher-order arctangent numbers. Then, the higher-order tangent numbers $\mathscr{T}(n,k)$ are given by}
\begin{equation}\mathscr{T}(n,k) = \frac{1}{k!}\cdot\frac{1}{(k-1)!} \sum_{r\,= 0}^{k-1}  T_{n + r} \,\mathscr{T}^{*}(k, r + 1)\quad(n,k\geq1).
\end{equation}

\vskip 2 mm \noindent {\em Demonstration.} Clearly, by the definition (1.1), we have that

\begin{equation*}\mathscr{T}(n,k)= \frac{1}{k!}\left. D_{t}^n \tan^{k} t \right|_{t = 0}^{}.
\end{equation*}

\vskip 2mm Now, since $T_n = \left. D_{t}^n \tan t \right|_{t = 0}^{}$, the proposed formula for $\mathscr{T}(n,k)$ in (2.15) could be obtained  upon setting $n = k$ in (2.14), then differentiating $n$ times both sides of the resulting expression with respect to $t$ and putting $t = 0$.

\vskip 2mm  In view of Proposition 7, the higher-order tangent numbers $\mathscr{T}(n,k)$ are linear combinations of the tangent numbers $T_n$ as is demonstrated through examples in Example 2.

\vskip 2mm \noindent{\bf Example 2.} The  numbers $\mathscr{T}(n,k)$  expressed through the tangent numbers $T_n$

\begin{small}
\begin{align*}&\mathscr{T}(1,1) = T_1,\quad \mathscr{T}(2,2) = \tfrac{1}{2}\,T_3\quad \mathscr{T}(3,1) = T_3,\quad\mathscr{T}(3,3) = -\tfrac{1}{6} \,T_3 + \tfrac{1}{12}\,T_5,
\\
&\mathscr{T}(4,2) = \tfrac{1}{2}\, T_5,\quad\mathscr{T}(4,4) = -\tfrac{1}{18} \,T_5 + \tfrac{1}{144}\, T_7,
\\
&\mathscr{T}(5,1) = T_5,\quad\mathscr{T}(5,3) =  - \tfrac{1}{6}\, T_5 + \tfrac{1}{12}\, T_7,\quad\mathscr{T}(5,5) = \tfrac{1}{120} T_5 - \tfrac{1}{144}\, T_7 + \tfrac{1}{2880} T_9,
\\
&\mathscr{T}(6,2) = \tfrac{1}{2} \,T_7,\quad\mathscr{T}(6,4) = -\tfrac{1}{18}\, T_7 + \tfrac{1}{144} \, T_9,\quad\mathscr{T}(6,6) = \tfrac{23}{10800}\,T_{7} - \tfrac{1}{2160}\,T_9 + \tfrac{1}{86400}\,T_{11},
\\
&\mathscr{T}(7,1) = T_7,\quad\mathscr{T}(7,3) = - \tfrac{1}{6} \,T_7 + \tfrac{1}{12} \,T_9,
\quad\mathscr{T}(7,5) = \tfrac{1}{120}\, T_7 - \tfrac{1}{144}\, T_9 + \tfrac{1}{2880}\,T_{11},
\\
& \mathscr{T}(7,7) = - \tfrac{1}{5040}\,T_7 + \tfrac{7}{32400}\,T_9 - \tfrac{1}{51840}\,T_{11} + \tfrac{1}{3628800} \,T_{13},
\\
&\mathscr{T}(8,2) = \tfrac{1}{2}\,T_9, \quad\mathscr{T}(8,4) = -\tfrac{1}{18}\, T_9 + \tfrac{1}{144} \, T_{11},
\quad\mathscr{T}(8,6) = \tfrac{23}{10800}\,T_{9} - \tfrac{1}{2160}\,T_{11} + \tfrac{1}{86400}\,T_{13},
\\
& \mathscr{T}(8,8)= - \tfrac{11}{264600}\,T_9+\tfrac{11}{907200}\,T_{11}-\tfrac{1}{1814400}\,T_{13}+\tfrac{1}{203212800}\,T_{15}.
 \end{align*}
\end{small}

\noindent{\bf Proposition 8.} {\em We have:}
\begin{small}
\begin{equation}
\mathscr{T}(n,k)= (-1)^{\frac{n-k}{2}} (-1)^n  \frac{2^n}{k!} \sum_{\alpha\,= k}^{n}\sum_{\beta\,=1}^{\alpha} (-1)^{\beta} \binom{\alpha-1}{k-1}\binom{\alpha}{\beta}\frac{\beta^n}{2^{\alpha}}\qquad(n\geq1,k\geq0).
\end{equation}
\end{small}

\vskip -2 mm \noindent {\em Demonstration.} To prove this formula, it suffices to recall the following power series expansion of $\tan^{k} x$ which was deduced by Schwatt \cite[p. 67, Eq. (76)]{Schwatt}
\begin{equation*}\tan^{k} x = \sum_{n\,= k}^{\infty} (-1)^{\frac{n-k}{2}} (-1)^n  \,2^n \,\sum_{\alpha\,= k}^{n}\sum_{\beta\,= 1}^{\alpha} (-1)^{\beta} \binom{\alpha-1}{k-1}\binom{\alpha}{\beta}\frac{\beta^n}{2^{\alpha}}\cdot \frac{x^n}{n!}.
\end{equation*}

\vskip 2mm \noindent{\bf Proposition 9.} {\em Let respectively $S(n,k)$ and $L_{n,k}$ be the Stirling numbers of the second kind and Lah numbers defined by} \cite[p. 50]{Comtet}
\begin{equation*}
\left(e^t-1\right)^k = k!\, \sum_{n\, = k}^{\infty}
S(n,k)\,\frac{t^n}{n!}
\end{equation*}
\noindent{\em and} \cite[p. 156]{Comtet}
\begin{equation*}
L_{n,k} = (-1)^n \binom{n-1}{k-1}\,\frac{n!}{k!}.
\end{equation*}

\vskip 2mm{\em Then, for $n\geq1$ and $k\geq0,$ we have that}
\begin{equation}
\mathscr{T}(n,k)= (-1)^{\frac{n-k}{2}} (-1)^n \sum_{\alpha\,= k}^{n} (-1)^{\alpha}\,2^{n -\alpha} \,S(n,\alpha)  \binom{\alpha-1}{k-1} \frac{\alpha!}{k!}
\end{equation}
\noindent{\em and}
\begin{equation}
\mathscr{T}(n,k)= (-1)^{\frac{n-k}{2}} (-1)^n \sum_{\alpha\,= k}^{n}  \,2^{n -\alpha} L_{\alpha,k}\, S(n,\alpha).
\end{equation}
\vskip 2 mm \noindent {\em Demonstration.} It is well--known that the Stirling numbers of the second kind are given by means of $S(n,k) = \frac{1}{k!} \sum_{j\,= 1}^k (-1)^{k-j} \binom{k}{j}\,j^n$ \cite[p. 204]{Comtet}. It is clear, then, that the sought formula (2.17) follows straightforwardly  by making use of this sum and (2.16). Furthermore, the expression on the right-side of the equation (2.17), in view of the definition of the Lah numbers $L_{n,k},$  may be written in the form given by (2.18).

\vskip 2mm For the sake of completeness of this paper, the following result  is  reproduced from the work of Butzer {\em et al.} \cite[Proposition 7.5, p. 482]{Butzer}, to which the interested reader is referred for a simple proof.

\vskip 2mm \noindent{\bf Proposition 10.} {\em In terms of central factorial numbers  $T(n,k)$ which are explicitly given by $T(n,k)=\frac{1}{k!}\sum_{\alpha\,= 0}^k (-1)^{\alpha} \binom{k}{\alpha}\left(\frac{k}{2}-\alpha\right)^n$ \textup{(}see \cite[p. 429]{Butzer} and \cite[pp. 213--217]{Riordan}\textup{)}, for $n,k\geq 1,$  we have:}
\begin{small}
\begin{align}
&\mathscr{T}(2 n,2 k) =  \sum_{\alpha\,= k}^{n} (-1)^{n - \alpha}\,2^{2 n -2 \alpha} \frac{(2 \alpha)!}{(2 k)!} \binom{\alpha -1}{\alpha - k} \,T(2 n,2 \alpha),
\\
&\mathscr{T}(2 n + 1,2 k +1) =  \sum_{\alpha\,= k}^{n} (-1)^{n - \alpha}\,2^{2 n -2 \alpha} \frac{(2 \alpha +1)!}{(2 k +1)!} \binom{\alpha -1/2}{\alpha - k} \,T(2 n +1,2 \alpha +1).
\end{align}
\end{small}

\section{Concluding remarks}

In this paper, two sequences of non-negative integer numbers, $\left\{\mathscr{T}(n,k)\right\}_{n,k\,=0}^{\infty}$ and $\left\{\mathscr{S}(n,k)\right\}_{n,k\,=0}^{\infty}$, have been thoroughly investigated and many of their properties are determined. In particular, it should be mentioned that the higher-order tangent numbers $\mathscr{T}(n,k)$ constitute a special class  of the partial multivariate Bell polynomials. Furthermore,  from the knowledge of  $\mathscr{T}(n,k)$, the higher-order secant numbers $\mathscr{S}(n,k)$ can be easily computed since $\mathscr{S}(n,k)$ and $\mathscr{T}(n,k)$ are related (see Proposition 3). In addition,  a simple explicit formula involving a double finite sum is deduced for the higher-order tangent numbers $\mathscr{T}(n,k)$ (Proposition 8) and it is shown that  $\mathscr{T}(n,k)$ are linear combinations of the classical tangent numbers $T_n$ (Proposition 7).

As an example of an application of the numbers $\mathscr{T}(n,k)$  and $\mathscr{S}(n,k)$ we give, without proof, recently established \cite[Corollaries 1 and 2]{Cvijovic}, and here slightly modified (observe that  $T(n,k) = \mathscr{T}(n,k) \,k!$  and  $S(n,k) = \mathscr{S}(n,k)\,k!$), closed-form higher derivative formulae.

\vskip 5mm \noindent{\bf Proposition 11.} {\em In terms of the tangent and  secant numbers of order} $k,$  $\mathscr{T}(n,k)$ {\em and} $\mathscr{S}(n,k),$ {\em for} $n\geq0,$ {{\em we have:}
\begin{align*}
& \qquad \textup{(a)}\quad D_{t}^{n}\,\tan x = \mathscr{T}(n,1)+ \sum_{k\,= 1}^{n + 1}(k-1)! \mathscr{T}(n +1,k)\tan^k x;\hskip 50mm
\\
&\qquad\textup{(b)}\quad D_{t}^{n}\,\sec x= \sec x\,\sum_{k\,=0}^{n} k!\, \mathscr{S}(n,k) \,\tan^k x;
\\
&\qquad\textup{(c)}\quad D_{t}^{n}\,\cot x= (-1)^n\,\left[\mathscr{T}(n,1) + \sum_{k\,= 1}^{n+1} (k-1)! \,\mathscr{T}(n+1,k)\,\cot^k x\right];
\\
&\qquad\textup{(d)}\quad D_{t}^{n}\,\csc x = (-1)^n\,\csc x \,\sum_{k\,= 0}^{n} k!\, \mathscr{S}(n,k)\,\cot^k x.
\end{align*}

\vskip 2mm \noindent{\bf Proposition 12.} {\em In terms of the tangent and  secant numbers of order} $k,$  $\mathscr{T}(n,k)$ {\em and} $\mathscr{S}(n,k),$ {\em for} $n\geq0,$ {{\em we have:}
\begin{align*}
& \qquad\textup{(a)}\quad  D_{t}^{n}\,\tanh x = (-1)^{\frac{n-1}{2}}\,\mathscr{T}(n,1) + \sum_{k\,= 1}^{n+1} (-1)^{\frac{n+k-1}{2}} \,(k-1)! \mathscr{T}(n +1,k)\tanh^k x;\hskip 50mm
\\
&\qquad\textup{(b)}\quad D_{t}^{n}\,\textup{sech\,} x = \textup{sech\,} x\,\sum_{k\,= 0}^{n}(-1)^{\frac{n+k}{2}}\, \mathscr{S}(n,k)\,\tanh^k x;
\\
&\qquad\textup{(c)}\quad  D_{t}^{n}\,\coth x = (-1)^{\frac{n-1}{2}}\,\mathscr{T}(n,1) + \sum_{k\,=1}^{n+1} (-1)^{\frac{n+k-1}{2}} \,(k-1)! \mathscr{T}(n +1,k)\,\coth^k x;
\\
&\qquad\textup{(d)}\quad  D_{t}^{n}\,\textup{csch\,}x = \textup{csch\,} x\,\sum_{k\,=0}^{n}(-1)^{\frac{n+k}{2}}\, \mathscr{S}(n,k)\,\coth^k x.
\end{align*}

We remark that, apart from $\mathscr{T}(n,k)$, which were introduced by Comtet (see the definition of $T(n,k)$ in \cite[p. 259]{Comtet}),  in the literature can be found several more definitions of higher-order tangent numbers: Carlitz and Scoville  considered  the numbers $T_{n}^{(k)}:= k! \mathscr{T}(n,k)$ (\cite[p. 428]{Carlitz} and \cite[p. 305]{Carlitz2}), while Lomont  studied  $C_{n}^{(k)}$ given by  $C_{n}^{(k)}:= (-1)^{\frac{n-k}{2}} (-1)^n \,k!\, \mathscr{T} (n,k)$ \cite[Eq. (3.1)]{Lomont}. However, it seems that  numbers $T_{n}^{(z)}$ defined by Cenkci \cite[p. 1500]{Cenkci} cannot be brought into connection with $\mathscr{T}(n,k)$. We note also that  Carlitz and Scoville  \cite[p. 428]{Carlitz} defined higher-order secant numbers  $S_{n}^{(k)}:= k! \mathscr{S}(n,k)$ (see also \cite[p. 665]{Knuth}). The tangent and secant numbers, $T_n$ and $S_n$, were first studied  by Andr\'e in 1879 \cite{Andre}.

The majority of results presented here are  new, while the known results are two recurrence relations given by Proposition 5 [for the recurrence (2.8) see \cite[Eq. (7)]{Knuth} and \cite[p. 258]{Comtet}; for the recurrence (2.9) see \cite[p. 665, Eq. (9)]{Knuth}] as well as the formula given by Proposition 6(b) [{\em cf.} (2.11) and \cite[p. 665, Eq. (10)]{Knuth}]. Also, the formula (2.18) was deduced earlier by different arguments \cite[p. 156]{Lomont2}. In addition, it should be noted that, judging by the results of extensive numerical calculations which support such a conclusion, our simple and rather compact  formulae (2.6) and (2.15) appear to be equivalent to the (relatively complicated-looking)  formulae which were established by Lomont \cite[see, respectively, Eq. (3.2) and Eqs. (3.4), (3.6) and (3.7)]{Lomont}. \\[-8mm]\section*{Acknowledgements}\noindent {\small The author acknowledges financial support from Ministry of Science of the Republic of Serbia under Research Projects 45005 and 172015.}

\end{document}